\documentclass[11pt,a4paper]{article}
\usepackage[T1]{fontenc}
\usepackage[utf8]{inputenc}
\usepackage{lmodern}
\usepackage[margin=25mm]{geometry}
\usepackage{amsmath,amssymb,amsthm,mathtools}
\usepackage{enumitem}
\usepackage{microtype}

\usepackage{xcolor}
\definecolor{revisionblue}{RGB}{0,70,170}
\newif\ifshowrevisions
\showrevisionstrue

\usepackage[hidelinks]{hyperref}
\newtheorem{theorem}{Theorem}[section]
\newtheorem{proposition}[theorem]{Proposition}
\newtheorem{lemma}[theorem]{Lemma}
\newtheorem{corollary}[theorem]{Corollary}
\theoremstyle{definition}

\newtheorem{example}[theorem]{Example}
\theoremstyle{remark}
\newtheorem{remark}[theorem]{Remark}
\DeclareMathOperator{\Aut}{Aut}
\DeclareMathOperator{\Sym}{Sym}
\newcommand{\id}{\mathrm{id}}

\newcommand{\cB}{\mathcal B}
\newcommand{\cH}{\mathcal H}

\newcommand{\doi}[1]{\href{https://doi.org/#1}{doi:\nolinkurl{#1}}}
\hypersetup{pdftitle={Stabilizer transfer and list colorings breaking graph automorphisms},
pdfauthor={Hui Lei, Yu Wang and Yueyue Zhao}}

\title{Colorings breaking prescribed families of graph automorphisms}
\author{Hui Lei\thanks{School of Statistics and Data Science, LPMC and KLMDASR,
Nankai University.
\texttt{hlei@nankai.edu.cn}.
Funded by the National Natural Science Foundation of China (Nos.~12431013), the Natural Science Foundation of Tianjin (24JCYBJC01670), and the Fundamental Research Funds for the Central Universities, Nankai University.}
\quad Yu Wang\thanks{Center for Combinatorics and LPMC, Nankai University.
\texttt{wangyu97@mail.nankai.edu.cn}.}
\quad Yueyue Zhao\thanks{School of Mathematical Sciences and LPMC, Nankai University and Department of Applied Mathematics, The Hong Kong Polytechnic University.
\texttt{799415027@qq.com}.}}
\date{}

\begin{document}

\maketitle

\begin{abstract}
Let \(\Gamma\) act on a set \(X\), and let \(\mathcal A\subseteq\Gamma\setminus\{\id\}\) be conjugacy invariant.
A coloring of \(X\) is \(\mathcal A\)-breaking if no element of \(\mathcal A\) preserves it. The ordinary parameter $D_{\mathcal A}(\Gamma,X)$ is the least positive integer $k$ admitting such a coloring $X\to[k]$; the list parameter $D_{\ell,\mathcal A}(\Gamma,X)$ is the least positive integer $k$ for which every assignment of $k$-element lists admits one.
For \(\mathcal A=\Gamma\setminus\{\id\}\), an
\(\mathcal A\)-breaking coloring is called a distinguishing coloring,
and we write \(D(G)\) and \(D_\ell(G)\) for the ordinary and list
parameters on \(V(G)\), respectively, and \(D'(G)\) and \(D'_\ell(G)\)
for the corresponding parameters on \(E(G)\).
We prove that for every continuous action of a compact group on a discrete set, each finite-color template admits a coloring from arbitrary lists of the same size whose stabilizer is contained in a conjugate of the template stabilizer.
Consequently, the ordinary and list parameters coincide whenever the ordinary parameter is finite.
In particular, we answer the question of Ferrara, Flesch and Gethner negatively for finite graphs by proving \(D_{\ell}(G)=D(G)\), and establish the finite case of the conjecture of Kwa\'sny and Stawiski by proving \(D'_{\ell}(G)=D'(G)\) whenever a distinguishing edge coloring exists.
We also show that every closed locally compact permutation group with infinite motion on a countably infinite set admits continuum many inequivalent distinguishing colorings from arbitrary 2-element lists.
Using compact kernels and finite bases, we extend ordinary and list equality to countable graphs with countable automorphism groups and prove that one additional list color suffices for connected locally finite graphs.
Using this multiplicity and a construction fixing all vertices of infinite degree, we prove that \(D'_\ell(G)\le 2\) for every countable connected graph with infinite edge motion, settling the list infinite edge-motion conjecture of Kwa\'sny and Stawiski.
An automorphism of a graph is small if it maps some vertex to a neighbor.
Finally, every connected countable graph outside \(\{K_2,K_3,K_4,K_5,C_5\}\) admits an edge coloring from arbitrary 2-element lists that breaks all small automorphisms, settling the conjecture of Kwa\'sny and Stawiski.
The latter conclusion extends to graphs in which each vertex has at most countably many neighbors in its own automorphism orbit.
\end{abstract}


\section{Introduction}\label{sec:introduction}
All graphs in this paper are simple, colorings need not be proper, and \emph{countable} means at most countable.
For an action of $\Gamma$ on $X$ (\(\Gamma\curvearrowright X\) ), put \( \Gamma_c=\{g\in\Gamma:c(gx)=c(x)\text{ for all }x\in X\}. \) Given a conjugacy-invariant family $\mathcal A\subseteq\Gamma\setminus\{\id\}$, a coloring is \emph{$\mathcal A$-breaking} if $\Gamma_c\cap\mathcal A=\emptyset$.
The ordinary parameter $D_{\mathcal A}(\Gamma,X)$ is the least positive integer $k$ admitting such a coloring $X\to[k]$; the list parameter $D_{\ell,\mathcal A}(\Gamma,X)$ is the least positive integer $k$ for which every assignment of $k$-element lists admits one.
Identical lists give $D_{\mathcal A}(\Gamma,X)\leq D_{\ell,\mathcal A}(\Gamma,X)$, and both parameters are $1$ when $\mathcal A=\emptyset$.
For \(\mathcal A=\Gamma\setminus\{\mathrm{id}\}\), an
\(\mathcal A\)-breaking coloring is called a \emph{distinguishing coloring},
and we write \(D(G)\) and \(D_\ell(G)\) for the ordinary and list
parameters on \(V(G)\), respectively, and \(D'(G)\) and \(D'_\ell(G)\)
for the corresponding parameters on \(E(G)\).
Write $[k]=\{1,\ldots,k\}$.

We first establish a stabilizer-transfer theorem for compact group actions, yielding equality between the ordinary and list parameters.

\begin{theorem}[Compact stabilizer transfer]\label{thm:compact}
Let a compact group $\Gamma$ act continuously on a discrete set $X$.
For every finite $k\ge1$, every template $f:X\to[k]$, and every assignment of $k$-element lists, there exist a legal coloring $c$ and an $h\in\Gamma$ such that
		\[
		\Gamma_c\le h\Gamma_fh^{-1}.
		\]
Consequently, for every conjugacy-invariant forbidden family $\mathcal A$,
		\[
		D_{\ell,\mathcal A}(\Gamma,X)=D_{\mathcal A}(\Gamma,X)
		\]
whenever the ordinary parameter is finite.
Neither faithfulness nor countability of $X$ is required.
	\end{theorem}

In particular, for every finite graph \(G\),
\[
 D_\ell(G)=D(G),\qquad D'_\ell(G)=D'(G),
\]
where the edge equality holds whenever a distinguishing edge coloring exists.

The proof maximizes discrete convex potentials over the compact template orbit.
The finite case answers the question of Ferrara, Flesch and Gethner~\cite[Question~1]{FFG} negatively for finite graphs and proves the finite list distinguishing index conjecture of Kwa\'sny and Stawiski~\cite[Conjecture~1.1]{KSindex}.
 The finite-orbit and rooted consequences are collected in Section~\ref{sec:graphs}.

The action has \emph{infinite motion} if every nonidentity element moves infinitely many points.
Lehner~\cite[Theorem~2]{LehnerLC} proved the ordinary two-color existence theorem for closed locally compact permutation groups with infinite motion on a countable set.
Our list version also provides continuum many inequivalent colorings.

\begin{theorem}[List infinite-motion theorem]
		\label{thm:list-infinite-motion}
Let $X$ be countable and $\Gamma\le\Sym(X)$ closed and locally compact in the permutation topology.
If $\Gamma$ has infinite motion, then every assignment of 2-element lists admits a distinguishing coloring.
If $X$ is infinite, there are $2^{\aleph_0}$ pairwise $\Gamma$-inequivalent such colorings.
Here $c$ and $d$ are equivalent if there is \(g\in\Gamma\) such that \(d(x)=c(g^{-1}x)\) for every \(x\in X\); the lists need not be $\Gamma$-invariant.
	\end{theorem}

Beyond compact actions, we separate a compact kernel from an invariant set whose induced action is distinguishable from 2-element lists.
This gives a stabilizer inclusion and an exact parameter formula (Theorem~\ref{thm:compact-kernel-motion}).
For countable finite-base actions, it also gives stabilizer transfer and ordinary/list equality (Theorem~\ref{thm:finitebase}).
The graph consequences include equality for countable graphs with countable automorphism groups and a bound of one additional list color for connected locally finite graphs; see Section~\ref{sec:graphs}.

Combining this multiplicity with a construction fixing all vertices of infinite degree yields the following list extension of Lehner's infinite edge-motion theorem~\cite{LehnerEdge}.

\begin{theorem}[Countable graphs with infinite edge motion]
		\label{thm:countable-edge-motion}
Let $G$ be a countable connected graph such that every nonidentity automorphism moves infinitely many edges.
Then every assignment of 2-element edge lists admits a distinguishing coloring.
In particular,
		\[D'_\ell(G)\le 2.
		\]
	\end{theorem}

Theorem~\ref{thm:countable-edge-motion} settles the conjecture of Kwa\'sny and Stawiski~\cite[Conjecture~8]{KSregular}.
The corresponding two-list conclusion for vertex colorings of connected locally finite graphs with infinite vertex motion follows from Corollary~\ref{cor:list-infinite-motion-graphs}.

Our final main result concerns \emph{small automorphisms}: following Kalinowski, Pil\'sniak and Wo\'zniak~\cite{KPW}, an automorphism $g$ is \emph{small} if $xg(x)\in E(G)$ for some vertex $x$.
Write $D_s,D_{\ell,s}$ for the vertex parameters and $D'_s,D'_{\ell,s}$ for the edge parameters.
Kwa\'sny and Stawiski~\cite{KSsmall} proved a 3-list edge bound for locally finite graphs without a $K_2$ component.
Stawiski~\cite{Stawiski} proved the ordinary 2-color bound, outside finite exceptions, for graphs of arbitrary cardinality.
We prove the 2-list conclusion under a condition weaker than local finiteness, settling the conjecture in~\cite[Conjecture~1]{KSsmall}.
Note that the graph \(C_4\) is not an exception to the small 2-list conclusion.

\begin{theorem}[2-element edge lists]\label{thm:smalltwo}
Let $G$ be a connected graph such that
		\begin{equation}\label{eq:countable-orbit-degree}
			|N_G(v)\cap\Aut(G)v|\le\aleph_0\qquad(v\in V(G)).
		\end{equation}
If $G\not\cong K_2,K_3,K_4,K_5,C_5$, every assignment of 2-element edge lists admits a coloring breaking all small automorphisms.
In particular, this holds for every connected countable graph outside the five exceptions, including graphs with infinite degrees.
	\end{theorem}

\section{Compact stabilizer transfer}\label{sec:finite}\label{sec:compact}
For a left action of $\Gamma$ on $X$, put $(g\cdot f)(x)=f(g^{-1}x)$, so that $\Gamma_{g\cdot f}=g\Gamma_fg^{-1}$.
Lists need not be invariant under $\Gamma$.

	\begin{example}[Full and small distinction]\label{gen:graph-families}
For $\Gamma=\Aut(G)$, the full and small forbidden families are
		\[
		\mathcal A_{\rm all}=\Gamma\setminus\{\id\},\qquad
		\mathcal A_{\rm sm}=\{g:xg(x)\in E(G)\text{ for some }x\}.
		\]
Their ordinary/list parameters are respectively
		\[
		\begin{array}{c|cc}
			&\mathcal A_{\rm all}&\mathcal A_{\rm sm}\\ \hline
			V(G)&(D,D_\ell)&(D_s,D_{\ell,s})\\
			E(G)&(D',D'_\ell)&(D'_s,D'_{\ell,s}).
		\end{array}
		\]
Both families are conjugacy invariant and open in the topology of pointwise convergence.
For the small family, conjugation carries an adjacent witness $x,g(x)$ to $hx,hg(x)$, and
		\[
		\mathcal A_{\rm sm}
		=\bigcup_x\ \bigcup_{y\in N_G(x)}\{g:g(x)=y\}
		\]
is a union of open point transporters.
	\end{example}

	\begin{lemma}[Passage to the permutation image]\label{gen:action-kernel}
Let $N$ be the kernel of the action and $\pi:\Gamma\to\Gamma|_X$ its permutation image.
If $\mathcal A\cap N\ne\emptyset$, no $\mathcal A$-breaking coloring exists.
Otherwise replacing $(\Gamma,\mathcal A)$ by $(\Gamma|_X,\pi(\mathcal A))$ preserves both parameters.
	\end{lemma}
	\begin{proof}
For every coloring $c$, its stabilizer is the inverse image under $\pi$ of its stabilizer in $\Gamma|_X$.
	\end{proof}

We use the standard compactness criterion for permutation groups on an arbitrary set; see Tornier~\cite[Section~2.1.1]{Tornier}.
For the countable setting, see also Lehner~\cite[Section~2]{LehnerLC}.
 Permutation groups carry the topology of pointwise convergence.

\begin{lemma}[\cite{{LehnerLC}, Tornier}]\label{lem:compactorbits}
A closed subgroup $\Gamma\leq \Sym(X)$ is compact if and only if all its point orbits are finite.
	\end{lemma}

\begin{lemma}[A discrete convexity inequality]\label{lem:convex}
Let \(\lambda_1<\cdots<\lambda_k\) be real numbers and put \(F(r)=\sum_{j=1}^r\lambda_j\) for \(1\leq r\leq k\).
For \(r,s\in[k]\),
		\[
		F(s)-F(r)\geq(s-r)\lambda_r,
		\]
and the inequality is strict if \(s>r\).
	\end{lemma}
	\begin{proof}
If \(s>r\), sum \(\lambda_{r+1},\ldots,\lambda_s\), each of which is larger than \(\lambda_r\).
If \(s<r\), then
		\[
		F(s)-F(r)=-\sum_{j=s+1}^r\lambda_j\geq-(r-s)\lambda_r.
		\]
The case \(s=r\) is equality.
	\end{proof}

	\begin{proof}[Proof of Theorem~\ref{thm:compact}]
Each point orbit is finite, being a compact discrete image of $\Gamma$.
The template orbit
		\[
		\mathcal K=\Gamma\cdot f
		\]
is compact: the map $g\mapsto g\cdot f$ is continuous into $[k]^X$.
In particular, $\mathcal K$ is the actual orbit, not merely its closure.

On each point orbit $O$, encode the finite palette $\bigcup_{x\in O}L(x)$ by an injection $\iota_O$ into $\mathbb R$.
Write $L(x)=\{a_{x,1},\ldots,a_{x,k}\}$ in increasing encoded order, with $\lambda_{x,r}=\iota_O(a_{x,r})$, and put $F_x(r)=\sum_{j\le r}\lambda_{x,j}$.
Well-order the point orbits as $(O_\alpha)_{\alpha<\theta}$ and set
		\[
		P_\alpha(t)=\sum_{x\in O_\alpha}F_x(t(x)).
		\]
Starting with $\mathcal K_0=\mathcal K$, let $\mathcal K_{\alpha+1}$ be the maximizers of $P_\alpha$ on $\mathcal K_\alpha$, and take intersections at limit stages.
Each potential is continuous with finite range.
Thus every stage is nonempty and compact, by attainment of the maximum at successor stages and compactness at limit stages.
Choose $t\in\mathcal K_\theta$ and set $c(x)=a_{x,t(x)}$.

Take $g\in\Gamma_c$.
On each finite point orbit $O$, color preservation and reindexing give
		\begin{equation*}\label{eq:orbitzero}
			\sum_{x\in O}\bigl((g\cdot t)(x)-t(x)\bigr)\lambda_{x,t(x)}=0.
		\end{equation*}
The rank sums of $g\cdot t$ and $t$ on $O$ also agree.
If the templates differ there, some rank increases, so Lemma~\ref{lem:convex} gives a strict increase in the potential on $O$.
If $g\cdot t\ne t$, choose the first differing orbit $O_\alpha$.
Agreement on all earlier orbits ensures that $g\cdot t$ survived their maximizations, so $g\cdot t\in\mathcal K_\alpha$.
The strict increase contradicts the choice of $\mathcal K_{\alpha+1}$.
Hence $g\cdot t=t$.
As $t=h\cdot f$ for some $h\in\Gamma$, we obtain $\Gamma_c\le\Gamma_t=h\Gamma_fh^{-1}$.

If $f$ breaks $\mathcal A$, conjugacy invariance makes $h\Gamma_fh^{-1}$ disjoint from $\mathcal A$, so $c$ also breaks $\mathcal A$.
Identical lists give the reverse parameter inequality.
	\end{proof}

	\begin{corollary}[Finite transfer]\label{prop:finite}
		\label{cor:finitebreaking}
Every action on a finite set admits the stabilizer transfer of Theorem~\ref{thm:compact}, and hence ordinary/list equality for every conjugacy-invariant forbidden family whenever the ordinary parameter exists.
	\end{corollary}
	\begin{proof}
Pass to the finite permutation image, which is compact, and lift the stabilizer inclusion to $\Gamma$.
	\end{proof}

\section{ Locally compact groups with infinite motion}\label{sec:motion}

	\begin{proof}[ Proof of Theorem~\ref{thm:list-infinite-motion}]
We adapt the decomposition in Lehner~\cite[Lemmas~5--7 and the proof of Theorem~2]{LehnerLC}, retaining additional finite orbits to encode inequivalent colorings.
For the trivial group, the assertions follow by choosing independently at each point.
Otherwise $X$ is infinite.
Choose a finite nonempty $X_0$ with compact open setwise stabilizer $H$ and enumerate the finite $H$-orbits as $S_0,S_1,\ldots$.
The existence of $X_0$ is the standard local-compactness criterion recalled in~\cite[Section~2]{LehnerLC}.
In this proof, $\Gamma_Y$ and $\Gamma_{(Y)}$ denote the setwise and pointwise stabilizers of $Y$, respectively.
Put
		\[
		\Gamma_n=\bigcap_{i\ge n}\Gamma_{S_i}.
		\]
We have $H\le\Gamma_n\le\Gamma_{n+1}$.
Each $\Gamma_n$ is a closed group with finite orbits (it preserves each tail orbit and their finite complement), hence compact by Lemma~\ref{lem:compactorbits}.

By the block construction in Lehner~\cite[Lemma~6]{LehnerLC}, for each $n$ the tail $\bigcup_{j\ge n}S_j$ can be partitioned into consecutive finite blocks $U_1(n),U_2(n),\ldots$, each a union of $S_j$, with
		\begin{equation}\label{eq:motion-blocks}
				(\Gamma_n)_{(U_{i+1}(n))}\le(\Gamma_n)_{(U_i(n))}.
		\end{equation}
For every infinite set $I\subseteq\mathbb N$, Lehner's Lemma~7~\cite{LehnerLC} states that the action of $\Gamma_n$ on $\bigcup_{i\in I}U_i(n)$ is faithful and admits an ordinary distinguishing 2-coloring.
We use these two lemmas without reproof.
Applying Theorem~\ref{thm:compact} to this action of the compact group $\Gamma_n$ gives a distinguishing coloring from any 2-element lists on the same union.

For a finite set $Y$ of size $|X_0|$, define
		\[
		H(Y)=\{g\in\Gamma:gX_0=Y\},\qquad
		\mathcal X_n=\{gX_0:g\in\Gamma_n\}.
		\]
Lehner~\cite[Lemma~5]{LehnerLC} gives the transporter properties used below: all elements of a nonempty $H(Y)$ send each $S_i$ to the same set, and $Y\in\mathcal X_n$ implies $H(Y)\subseteq\Gamma_n$.
Enumerate the remaining finite sets $Y$ as $Y_1,Y_2,\ldots$, using dummy empty transporters if needed, and choose a sequence $n_k$ in which every nonnegative integer occurs infinitely often.

Reserve disjoint finite sets in stages.
At stage $k$, let $F_k$ be the union of earlier reservations.
If $H(Y_k)\ne\emptyset$, choose $g_k\in H(Y_k)$.
Since $g_k\notin\bigcup_n\Gamma_n$, infinitely many $i$ satisfy $g_kS_i\ne S_i$.
Choose such a $j_k$ with
		\[
		Z_k=S_{j_k}\cup g_kS_{j_k}
		\]
disjoint from $F_k$: only finitely many of the pairwise disjoint $S_i$ or their images meet $F_k$.
Otherwise put $Z_k=\emptyset$.
Next reserve one block $U_{i_k}(n_k)$ disjoint from $F_k\cup Z_k$.
Finally reserve an entire $H$-orbit $R_k$ disjoint from $F_k\cup Z_k\cup U_{i_k}(n_k)$.
Each of these choices is possible because the excluded set is finite and the candidate blocks or orbits form an infinite disjoint family.
Include all three sets among the reservations at subsequent stages.
Thus each $n$ receives infinitely many blocks, all reservations are pairwise disjoint, and the sets $R_k$ remain available for independent choices.

For each $n$, distinguish $\Gamma_n$ on the union of its reserved blocks using the preceding compact argument.
On each nonempty $Z_k$, choose colors making the color multisets on $S_{j_k}$ and $g_kS_{j_k}$ different.
To do this, leave one point of $S_{j_k}\setminus g_kS_{j_k}$ until last; at most one of its two colors can make the multisets equal.

On each $R_k$, fix the colors at all but one point and use the two choices at the remaining point to encode a bit.
These two colorings have different finite color multisets.
Fix all other choices, including the colors on points outside the reserved sets, independently of the bits.
We obtain a legal coloring $c_\varepsilon$ for every $\varepsilon\in\{0,1\}^{\mathbb N}$.

A nonidentity color preserver belongs to no $\Gamma_n$, because the reserved union for $n$ is $\Gamma_n$-invariant and distinguishes that group.
Its image of $X_0$ must therefore be some $Y_k$.
It maps $S_{j_k}$ onto $g_kS_{j_k}$, contradicting the different color multisets.
Thus every $c_\varepsilon$ is distinguishing, regardless of the choices on the $R_k$.

Distinct bit sequences give $H$-inequivalent colorings, since $H$ preserves each $R_k$ and hence its color multiset.
Moreover, $[\Gamma:H]\le\aleph_0$: a coset $gH$ is determined by the finite set $gX_0$.
Each $\Gamma$-orbit of colorings therefore contains at most countably many $H$-orbits.
The family just constructed consequently meets $2^{\aleph_0}$ distinct $\Gamma$-orbits, as required.
	\end{proof}

\section{ Compact kernels and finite-base actions}
	\label{sec:countableextension}
A \emph{base} for an action is a set whose pointwise stabilizer is trivial.
We write $D(\Gamma,X)$ and $D_\ell(\Gamma,X)$ for full distinction.

The countable infinite-motion lemma of Imrich, Smith, Tucker and Watkins~\cite[Lemma~3.3]{IMSTW} gives an ordinary distinguishing 2-coloring for a countably infinite permutation group with infinite motion on a countably infinite set.
The same successive-choice idea gives the following list version for a countable family of permutations; we include the short argument.

	\begin{lemma}[Countable list motion]\label{lem:listmotion}
Every countable family of infinite-support permutations of a set $W$ can be broken by a coloring from arbitrary 2-element lists.
	\end{lemma}
	\begin{proof}
Enumerate the family.
For each successive permutation $g$, choose $x,gx$ distinct and still uncolored; infinite support makes this possible, since only finitely many points have been colored.
Give these points different permissible colors.
Finally color all remaining points arbitrarily.
	\end{proof}

	\begin{theorem}[Compact kernels]
		\label{thm:compact-kernel-motion}\label{thm:compact-kernel-forbidden}\label{cor:locally-compact-kernel}
Let a Hausdorff topological group $\Gamma$ act continuously on a discrete set $X=W\sqcup U$, where $W$ is invariant and $K=\Gamma_{(W)}$ is compact.
Write $\pi:\Gamma\to\Gamma|_W$ for restriction.
		\begin{enumerate}[label=\textup{(\roman*)}]
\item If $\Gamma|_W$ is countable and has infinite motion, then for every finite $k\ge2$, template $f:X\to[k]$, and $k$-element lists, there are a legal coloring $c$ and $h\in K$ such that
			\[
			\Gamma_c\le h(K\cap\Gamma_f)h^{-1}.
			\]
\item Let $\mathcal A\subseteq\Gamma\setminus\{\id\}$ be nonempty and conjugacy invariant.
If $\pi(\mathcal A\setminus K)$ is countable and consists of infinite-support permutations, then
			\begin{equation}\label{eq:kernel-formula}
				D_{\ell,\mathcal A}(\Gamma,X)=D_{\mathcal A}(\Gamma,X)
				=\max\{2,D_{\mathcal A\cap K}(K,U)\},
			\end{equation}
whenever the right-hand parameter is finite; otherwise neither global parameter is finite.
\item If \(W\) is countable and \(\Gamma|_W\) is closed and locally compact with infinite motion, the stabilizer inclusion in \textup{(i)} holds.
For every nonempty conjugacy-invariant family \(\mathcal A\), formula~\eqref{eq:kernel-formula} also holds, with the same finiteness convention.
\end{enumerate}
	\end{theorem}
	\begin{proof}
For (i), Lemma~\ref{lem:listmotion} distinguishes $\Gamma|_W$ using 2-element sublists.
Compact transfer on $U$ gives $K_{c|_U}\le hK_{f|_U}h^{-1}$ for some $h\in K$.
Every global preserver lies in $K$, and $K_{f|_U}=K\cap\Gamma_f$.

For (ii), put $d=D_{\mathcal A\cap K}(K,U)$.
Given lists of size $\max\{2,d\}$, use Lemma~\ref{lem:listmotion} on $W$ to break $\pi(\mathcal A\setminus K)$ and compact transfer on $U$ to break $\mathcal A\cap K$.
Their union breaks $\mathcal A$.
Conversely, nonempty $\mathcal A$ requires at least two colors, and every global breaking coloring restricts to an $(\mathcal A\cap K)$-breaking coloring of $U$.
These lower bounds also cover nonexistence of a finite $d$; identical lists finish the finite case.

For \textup{(iii)}, use Theorem~\ref{thm:list-infinite-motion} to distinguish \(\Gamma|_W\) from two-element sublists, and repeat the compact-kernel argument in \textup{(i)}.
 Coloring \(U\) to break \(\mathcal A\cap K\), together with the same lower bounds as in \textup{(ii)}, gives~\eqref{eq:kernel-formula}.
\end{proof}

For the finite-base transfer theorem, we use the following criterion
to obtain infinite motion on the union of the infinite orbits.
\begin{lemma}[Finite bases and motion]\label{lem:basemotion}
\leavevmode
\begin{enumerate}[label=\textup{(\roman*)}]
\item Let \(Q\) act on \(W\), and let \(S,T\subseteq W\) be finite.
If every point orbit of \(Q\) on \(W\) has more than \(|S||T|\)
elements, then there exists \(q\in Q\) such that
\(qS\cap T=\varnothing\).
\item A faithful finite-base action with only infinite orbits has infinite motion.
\end{enumerate}
\end{lemma}
	
	\begin{proof} For (i), see ~\cite[Theorem~2]{BBOMN}. 
For (ii),  let $B$ be a finite base.
If a nonidentity element $g$ had finite support $S$, the cited separation lemma would give $q\in Q$ with $qS\cap B=\emptyset$.
The nonidentity conjugate $qgq^{-1}$ would then fix $B$ pointwise, contradicting the definition of a base.
	\end{proof}

	\begin{theorem}[Finite-base transfer]\label{thm:finitebase}
Let a countable group $\Gamma$ act faithfully on a nonempty set $X$ with a finite base.
For every conjugacy-invariant $\mathcal A\subseteq\Gamma\setminus\{\id\}$,
		\[
		D_{\ell,\mathcal A}(\Gamma,X)=D_{\mathcal A}(\Gamma,X)<\infty.
		\]
If $\Gamma$ is infinite, let $U$ be the union of its finite point orbits, $W=X\setminus U$, and $K=\Gamma_{(W)}$.
Then $W\ne\emptyset$, $K$ is finite, and for nonempty $\mathcal A$,
		\begin{equation}\label{eq:baseformula}
			D_{\ell,\mathcal A}(\Gamma,X)=D_{\mathcal A}(\Gamma,X)
			=\max\{2,D_{\mathcal A\cap K}(K,U)\}.
		\end{equation}
 For this action, every finite-color template \(f\) and lists of the same size admit a legal coloring \(c\) and \(h\in\Gamma\) with \(\Gamma_c\le h\Gamma_fh^{-1}\).
\end{theorem}
	\begin{proof}
A finite base has a finite distinguishing coloring: assign its points distinct colors and use one further color elsewhere.
Thus a finite ordinary parameter exists.
The empty family is immediate, and for finite $\Gamma$ compact transfer applies with the discrete topology.

Suppose $\Gamma$ is infinite and fix a finite base $B$.
Put $C=\Gamma(B\cap U)$ and $H=\Gamma_{(B\cap W)}$.
The finite invariant set $C$ gives an embedding $H\hookrightarrow\Sym(C)$, since an element of $H$ fixing $C$ fixes $B$.
Hence $H$ and $K\le H$ are finite, and $W$ cannot be empty.
The faithful image $Q=\Gamma|_W$ has infinite point orbits and a finite base: the stabilizer of $B\cap W$ is a finite image of $H$, so add a moved point for each of its nonidentity elements.
Lemma~\ref{lem:basemotion} gives infinite motion for $Q$.
The finite group $K$ acts faithfully on $U$, so its ordinary parameter is finite.
Theorem~\ref{thm:compact-kernel-motion}, with $\Gamma$ discrete, now gives~\eqref{eq:baseformula}.
 For the stabilizer inclusion, use compact transfer for finite \(\Gamma\) and Theorem~\ref{thm:compact-kernel-motion}(i) otherwise.
A one-color template has stabilizer \(\Gamma\), so that case is immediate.
\end{proof}

\section{ Graph applications}\label{sec:graphs}
Graph automorphism groups carry the topology of pointwise convergence on vertices.
The group $\Aut(G)$ is closed in $\Sym(V(G))$, because failure to preserve adjacency or nonadjacency is witnessed on two vertices.
The vertex and edge actions are continuous on discrete sets.

\subsection{ Finite orbits and one additional list color}

\begin{lemma}[Edge-action kernel]\label{lem:edgekernel}
The action of $\Aut(G)$ on $E(G)$ is faithful exactly when $G$ has no $K_2$ component and at most one isolated vertex.
Its kernel contains a small automorphism exactly when $G$ has a $K_2$ component.
	\end{lemma}

\begin{proof}
An automorphism fixing every edge fixes every vertex of degree at least two, as the intersection of two incident edges, and then every leaf outside a $K_2$ component.
The kernel therefore consists precisely of permutations of isolated vertices and independent endpoint swaps in $K_2$ components.
Only the latter can move a vertex to a neighbor.
	\end{proof}

\begin{corollary}[Finite orbits and rooted graphs]\label{cor:finiteorbits}
If every vertex orbit of $\Gamma=\Aut(G)$ is finite, then for $X=V(G)$ or $E(G)$ and every conjugacy-invariant forbidden family,
		\[
		D_{\ell,\mathcal A}(\Gamma,X)=D_{\mathcal A}(\Gamma,X)
		\]
whenever the ordinary parameter is finite.
The same holds for $\Gamma=\Aut(G,r)$ when $G$ is connected and locally finite.
	\end{corollary}

\begin{proof}[Proof of Corollary~\ref{cor:finiteorbits}]
Use Lemma~\ref{lem:compactorbits} and Theorem~\ref{thm:compact}.
In a connected locally finite rooted graph, root-stabilizer orbits lie in finite distance spheres.
	\end{proof}

\begin{corollary}\label{cor:finiteinvariant}
If a connected locally finite graph has a finite vertex orbit, then all its vertex orbits are finite, so Corollary~\ref{cor:finiteorbits} applies.
	\end{corollary}

\begin{proof}
The finite distance neighborhoods of that orbit are invariant and exhaust the vertex set.
	\end{proof}

\begin{proposition}[One marked point]\label{prop:marked}
Suppose $\Gamma$ acts continuously on discrete $X$ and the stabilizer of some $x_*\in X$ is compact.
Let $\mathcal A$ be conjugacy invariant.
If an $\mathcal A$-breaking $k$-coloring exists, then every $(k+1)$-list assignment admits an $\mathcal A$-breaking coloring.
	\end{proposition}

\begin{proof}
Choose $a\in L(x_*)$ and $k$-element sublists avoiding $a$ at all other points.
For $H=\Gamma_{x_*}$, the restricted template on $X\setminus\{x_*\}$ breaks $\mathcal A\cap H$.
Compact transfer colors that complement from its sublists, breaking $\mathcal A\cap H$.
Give $x_*$ color $a$; its unique color forces every preserver into $H$, completing the proof.
	\end{proof}

\begin{theorem}[One additional list color]\label{thm:plusone}
For a connected locally finite graph $G$, a conjugacy-invariant $\mathcal A\subseteq\Aut(G)\setminus\{\id\}$, and $X=V(G)$ or $E(G)$,
		\[
		D_{\mathcal A}(\Aut(G),X)\le D_{\ell,\mathcal A}(\Aut(G),X)
		\le D_{\mathcal A}(\Aut(G),X)+1
		\]
whenever the ordinary parameter is finite.
	\end{theorem}

\begin{proof}[Proof of Theorem~\ref{thm:plusone}]
In a connected locally finite graph, vertex stabilizers are compact.
An edge stabilizer is a finite extension of the compact pointwise stabilizer of its endpoints.
Apply Proposition~\ref{prop:marked}.
If the coordinate set is empty, ordinary existence forces $\mathcal A=\emptyset$, and both parameters are $1$.
	\end{proof}

\subsection{ Countable automorphism groups}
To apply Theorem~\ref{thm:finitebase} to both vertex and edge colorings, we first
establish the existence of finite bases for the two actions.

	\begin{lemma}[Countable closed groups and edge bases]
		\label{lem:closedbase}\label{lem:edgebase}
\leavevmode
\begin{enumerate}[label=\textup{(\roman*)}]
  \item A countable closed subgroup of $\Sym(X)$, with $X$ countable, has a finite base.
  \item If $\Aut(G)$ has a finite vertex base, then its faithful image on $E(G)$ also has a finite base. 
\end{enumerate}
 \end{lemma}   
	\begin{proof}
   For (i),  see \cite[Theorem~3.6(i)]{IMSTW}.
    For (ii), fix a finite vertex base $B$.
For each $v\in B$ of degree at least two, select two incident edges; for a leaf whose neighbor has degree at least two, select its edge and another edge at that neighbor.
Also select the edge of every $K_2$ component meeting $B$.
A permutation fixing these edges fixes all points of $B$ outside isolated vertices and $K_2$ components.
An element of the edge kernel can undo its action on the isolated vertices and the selected $K_2$ components (Lemma~\ref{lem:edgekernel}).
The resulting product fixes $B$ and is the identity.
Thus the selected edges form a base modulo the kernel.
	\end{proof}

	\begin{theorem}[Countable automorphism groups]\label{thm:countablegraphs}
Let $G$ be nonempty with both $V(G)$ and $\Gamma=\Aut(G)$ countable, and let $\mathcal A\subseteq\Gamma\setminus\{\id\}$ be conjugacy invariant.
Then
		\[
		D_{\ell,\mathcal A}(\Gamma,V(G))=D_{\mathcal A}(\Gamma,V(G))<\infty.
		\]
On edges, an $\mathcal A$-breaking coloring exists exactly when $\mathcal A\cap\ker(\Gamma\curvearrowright E(G))=\emptyset$; under this condition the two edge parameters are equal and finite.
Neither connectedness nor local finiteness is required.
	\end{theorem}

\begin{proof}[Proof of Theorem~\ref{thm:countablegraphs}]

Lemma~\ref{lem:closedbase} gives a finite vertex base, and Theorem~\ref{thm:finitebase} proves the vertex assertion.
For edges let $N$ be the action kernel.
Every coloring stabilizer contains $N$, so $\mathcal A\cap N=\emptyset$ is necessary.
Under this condition, pass to the countable faithful edge image and the image of $\mathcal A$.
Their parameters are unchanged, and Lemma~\ref{lem:edgebase} supplies a finite base.
Apply Theorem~\ref{thm:finitebase}; if $E(G)$ is empty, the condition forces $\mathcal A=\emptyset$ and both parameters are one.
	\end{proof}

In particular, countable graphs with countable automorphism groups satisfy $D_\ell=D$ and $D_{\ell,s}=D_s$.
The edge equalities hold for full distinction when the edge action is faithful, and for small distinction precisely when there is no $K_2$ component.

\subsection{ Locally finite graphs with infinite motion}

\begin{corollary}[Locally finite graphs with infinite motion]
		\label{cor:list-infinite-motion-graphs}
Let $G$ be connected, infinite, and locally finite.
If every nonidentity automorphism moves infinitely many vertices, arbitrary 2-element vertex lists and arbitrary 2-element edge lists admit distinguishing colorings.
For every nonempty conjugacy-invariant forbidden family, both ordinary and list parameters, on vertices or edges, equal two.
For a trivial automorphism group the full parameters are one.
	\end{corollary}

\begin{proof}
The vertex action is closed and has compact open vertex stabilizers.
The edge action is faithful, and its edge stabilizers are likewise compact and open.
Its image is closed: an edge-permutation limit can, after fixing the image of one edge, be taken inside a coset of its compact stabilizer.
Finally, every moved vertex is incident to a moved edge, since otherwise it is either the intersection of two fixed edges or belongs to a $K_2$ component.
Thus infinite vertex motion implies infinite edge motion.
Apply Theorem~\ref{thm:list-infinite-motion} to both actions; a nonempty forbidden family rules out one color.
	\end{proof}

\subsection{Countable graphs with infinite edge motion}
	
We now remove local finiteness from the edge-coloring conclusion.
For a partial edge coloring that colors every edge incident with a vertex $v$, its \emph{incidence color multiset} records, for each color, the number of edges of that color incident with $v$; these multiplicities may be infinite.

	\begin{lemma}[Fixing vertices of infinite degree]
		\label{lem:infinite-degree-multisets}
Let $G$ be a countable graph with 2-element lists on its edges, and put $I=\{v\in V(G):d_G(v)=\infty\}$.
The edges incident with $I$ admit a legal coloring such that distinct vertices of $I$ have distinct incidence color multisets.
Consequently, every color-preserving automorphism of any completion fixes $I$ pointwise.
	\end{lemma}
	\begin{proof}
Enumerate $I$ as $v_1,v_2,\ldots$, stopping if $I$ is finite.
At stage $n$, color all still uncolored edges incident with $v_n$, without changing any previous color.
At most $n-1$ incident edges are already colored, since each such edge joins $v_n$ to an earlier vertex of $I$.
Thus infinitely many incident edges remain uncolored.

These remaining lists permit infinitely many distinct incidence color multisets, even after accounting for the finitely many prescribed colors.
Indeed, if a color $a$ occurs in infinitely many remaining lists, use it on exactly $m$ remaining edges and avoid it on all others, for any $m\ge0$.
Its total multiplicity is a fixed finite number plus $m$.
Otherwise each color occurs in only finitely many remaining lists.
We can then choose infinitely many pairwise disjoint 2-element lists: each finite collection of previously chosen colors excludes only finitely many remaining edges.
Fix the colors on all other edges and vary the choices on the selected edges.
Different choices give different multisets, since the fixed contribution of every selected color is finite and that color occurs in no other selected list.

Choose a completion at $v_n$ whose multiset differs from those already assigned to $v_1,\ldots,v_{n-1}$.
Later stages do not change any edge incident with $v_n$.
The resulting multisets on $I$ are pairwise distinct.
Every automorphism preserves vertex degrees, and every color preserver also preserves these multisets, so it fixes $I$ pointwise.
	\end{proof}

	\begin{proof}[Proof of Theorem~\ref{thm:countable-edge-motion}]
Fix arbitrary 2-element lists on $E(G)$.
If $G$ is finite, the infinite-edge-motion assumption forces $\Aut(G)$ to be trivial, and any legal coloring suffices.
Hence suppose that $G$ is infinite.

Let $I$ be the set of vertices of infinite degree.
Apply Lemma~\ref{lem:infinite-degree-multisets} to color all edges incident with $I$.
Every color preserver of the completed coloring will therefore fix $I$ pointwise.
Each component $C$ of $G-I$ is locally finite, and each of its vertices has finite degree in the original graph $G$.
Label these vertices by
		\begin{equation}\label{eq:boundary-label}
			\beta(v)=N_G(v)\cap I.
		\end{equation}
The label records the actual neighbors in $I$, not just their number.
Every automorphism fixing $I$ induces label-preserving isomorphisms between components of $G-I$.

First consider a finite component $C$.
Any label-preserving automorphism of $C$ extends by the identity outside $C$ to an automorphism of $G$.
If it were nonidentity, it would move only finitely many vertices of finite degree in $G$, and hence only finitely many edges, a contradiction.
Moreover, two distinct finite components cannot be label-isomorphic: exchanging them by such an isomorphism and its inverse, and fixing the remainder of $G$, gives the same contradiction.
Thus every automorphism fixing $I$ fixes all finite components pointwise.

For an infinite component $C$, put
		\[
		\Lambda_C=\Aut(C,\beta).
		\]
This is a closed subgroup of $\Aut(C)$.
Every element of $\Lambda_C$ extends by the identity outside $C$ to an automorphism of $G$, because it preserves the labels in~\eqref{eq:boundary-label}.
We claim that every nonidentity element of $\Lambda_C$ moves infinitely many internal edges of $C$.
In an infinite connected graph, every moved vertex is incident with a moved edge: a vertex of degree at least two is fixed if two incident edges are fixed, and a leaf is fixed if its incident edge is fixed, since its neighbor has degree at least two and the endpoints cannot be interchanged.
Consequently, an element moving only finitely many internal edges moves only finitely many vertices.
All these vertices have finite degree in $G$, so its extension would move only finitely many edges of $G$, contrary to the hypothesis.

The edge action of $\Lambda_C$ is faithful, and its permutation image is closed and locally compact.
To verify these topological assertions, the setwise stabilizer of an internal edge is compact in the vertex permutation topology: its pointwise endpoint stabilizer is a closed subgroup of a compact vertex stabilizer, and has index at most two in the setwise edge stabilizer.
The edge action is continuous, so this stabilizer has compact image, which is open in the edge permutation topology.
For closedness, consider a convergent sequence of edge permutations induced by $\Lambda_C$.
The image of one fixed edge is eventually constant, so the corresponding automorphisms eventually lie in a fixed coset of its compact stabilizer.
A convergent subsequence in that coset induces the specified edge-permutation limit.
This proves closedness, since the permutation space on the countable edge set is metrizable.

By Theorem~\ref{thm:list-infinite-motion}, the lists on $E(C)$ admit $2^{\aleph_0}$ pairwise $\Lambda_C$-inequivalent distinguishing colorings.
We use this supply to prevent component exchanges.
Enumerate the infinite components as $C_1,C_2,\ldots$, again stopping if there are only finitely many, and color their internal edges in this order.
Require that each coloring $c_n$ distinguish $\Lambda_{C_n}$ and that the labeled edge-colored graphs $(C_n,\beta,c_n)$ be pairwise nonisomorphic.

To justify the choice at stage $n$, consider an earlier component $C_i$.
If $(C_n,\beta)$ and $(C_i,\beta)$ are not isomorphic, there is no restriction.
Otherwise choose a label-preserving isomorphism $\theta_i:C_n\to C_i$ and pull $c_i$ back to $C_n$.
All colorings of $C_n$ that make these two labeled edge-colored graphs isomorphic form exactly one $\Lambda_{C_n}$-orbit, namely the orbit of this pullback.
Hence the earlier components exclude at most $n-1$ orbits.
Choose a legal distinguishing coloring from any other orbit supplied by Theorem~\ref{thm:list-infinite-motion}.
No invariance of the lists under $\theta_i$ or $\Lambda_{C_n}$ is needed: an excluded orbit may contain no legal coloring, in which case it imposes no restriction.

Finally, color the internal edges of all finite components arbitrarily from their lists.
These assignments, together with the coloring of edges incident with $I$, use disjoint edge sets and cover $E(G)$.
Let $\varphi$ preserve the completed coloring.
It fixes $I$ pointwise by Lemma~\ref{lem:infinite-degree-multisets}, and therefore induces label-preserving isomorphisms between the components of $G-I$.
It fixes every finite component pointwise by the first part of the proof.
The pairwise nonisomorphism of the colored infinite components forces it to preserve each of them setwise; their internal distinguishing colorings then force pointwise fixation.
Thus $\varphi=\id$, proving $D'_\ell(G)\le2$.
	\end{proof}

	\begin{remark}
Under the assumptions of Theorem~\ref{thm:countable-edge-motion}, $D'_\ell(G)=D'(G)=2$ if $\Aut(G)$ is nontrivial, and both parameters are one otherwise.
The use of inequivalent colorings is essential to the component argument: distinguishing each infinite component internally does not by itself rule out exchanges of distinct components.
	\end{remark}

\section{Breaking small automorphisms with 2 edge-list colors}
\label{sec:regular}

To break small automorphisms, we will use the following result on
distinguishing edge colorings of regular graphs. Put $\cH=\{K_2,K_3,K_4,K_5,K_{3,3},C_4,C_5\}$. 
\begin{theorem}[\cite{KSregular}, Theorems~3 and~4]\label{thm1:regular}
Every connected countable regular graph of order at least seven  admits a distinguishing edge coloring from arbitrary 2-element lists.
\end{theorem}

\begin{lemma}\label{thm:regular}\label{lem:regularsmall}
Every connected countable regular graph outside $\cH$ admits a distinguishing edge coloring from arbitrary 2-element lists.
\end{lemma}
\begin{proof} 
Let \(G\) be a connected countable \(d\)-regular graph outside
\(\mathcal H\). If \(d=\aleph_0\), or if \(d\) is finite and \(G\)
has at least seven vertices, the conclusion follows from Theorem~\ref{thm1:regular}.
We may therefore assume that \(G\) has at most six vertices.

The case \(G=K_1\) is immediate, and the only connected \(1\)-regular
graph is \(K_2\), which belongs to \(\mathcal H\). For \(d=2,3\),
the conclusion follows from the base cases in the proof of
\textup{[9, Theorem~3]}. For \(d=5\), we have \(G=K_6\), and
\textup{[9, Lemma~10]} applies.

It remains to consider \(d=4\). Since \(K_5\in\mathcal H\), the graph
\(G\) has six vertices. Its complement is a perfect matching, so
\(G\cong K_{2,2,2}\). Label the three parts
\(\{a,y\}\), \(\{b,x\}\), and \(\{c,z\}\). Color the edges of the
triangle \(abc\), together with \(ax,by,yz\), red, and all other
edges blue. In the red subgraph, \(abc\) is the unique triangle.
Its vertices \(a,b,c\) are distinguished by their attached paths
of lengths \(1,2,0\), respectively. Hence every color-preserving
automorphism fixes the triangle pointwise and then fixes \(x,y,z\).
This gives a distinguishing two-edge-coloring of \(G\). By
Theorem~1.1, \(G\) also admits a distinguishing edge coloring from
arbitrary \(2\)-element lists.
\end{proof}

\begin{proposition}\label{prop:c4k33}
                The ordinary and list small edge parameters of $C_4$ and $K_{3,3}$
                are two.
        \end{proposition}
        \begin{proof}
                Color a two-edge path with one color and all remaining edges with
                another. Every preserving automorphism fixes its middle vertex, hence
                preserves the two bipartition classes and is not small. Apply finite
                transfer. Both graphs have a small automorphism, giving the lower bound.
        \end{proof}
	\begin{lemma}[A freely choosable edge]\label{lem:freeedge}\label{cor:avoid-colored-graph}
For $H\in\cH$, a vertex $r$, and 2-element edge lists, one can color all but an edge $e_*$ so that either choice on $e_*$ yields a coloring whose only preserving automorphism fixing $r$ is the identity.
The two completions have different numbers of $p$-colored edges, where $L(e_*)=\{p,q\}$.
 Consequently, for any fixed finite edge-colored graph \(F\), one can choose a completion for which no color-preserving isomorphism \(H\to F\) exists.
\end{lemma}
	\begin{proof}
Two edges can always receive different colors.
We also use this operation: between two individually preserved pairs, color one edge $p$ and avoid $p$ on the other three.
The unique $p$-edge fixes its endpoints and therefore both pairs pointwise.

For $K_2$, fixing $r$ suffices, so leave its edge free.
For $K_3,C_4,C_5$, give the two edges at $r$ different colors and leave any edge not incident with $r$ free.

For $K_4$ and $K_5$, color the edges at $r$ nonmonochromatically.
If all remaining vertices are fixed, leave any remaining edge free.
If the color classes leave a single pair \(\{a,b\}\), give \(va,vb\) different colors for a fixed vertex \(v\ne r\), and leave \(ab\) free.
For \(K_5\), the remaining class sizes are \(2+2\) and \(3+1\).
For $2+2$, apply the pair operation and leave an edge within a pair free.
For $3+1$, color the edges from the fixed vertex to the triple nonmonochromatically.
This either fixes the triple, or leaves a pair and a fixed vertex within it; separate that pair in the same way and leave its internal edge free.

For $K_{3,3}$ with parts $\{r,a_1,a_2\}$ and $\{b_1,b_2,b_3\}$, color the three edges at $r$ nonmonochromatically.
If the $b_i$ are fixed, give $a_1b_1,a_2b_1$ different colors and leave $a_1b_2$ free.
Otherwise assume $\{b_1,b_2\}$ is a color class.
Apply the pair operation to it and $\{a_1,a_2\}$, leaving $a_1b_3$ free.

Color other unused edges arbitrarily.
If $L(e_*)=\{p,q\}$, the two completions have different finite numbers of $p$-colored edges and are therefore nonisomorphic.
For the last assertion, at most one number equals the number of $p$-colored edges in $F$, so the other completion suffices.
\end{proof}

	\subsection{Orbit components and anchors}\label{sec:anchors}
Let $\Gamma=\Aut(G)$.
The connected components of $G[X]$, over all vertex orbits $X$ of $\Gamma$, are the \emph{orbit components}; write $\cB$ for this family and $X(H)$ for the original orbit containing $H$.
Let $E_{\mathrm{int}}=\bigcup_{H\in\cB}E(H)$, and let $E_G(H,J)$ denote the edges between distinct components.
Their quotient graph $Q$ has vertex set $\cB$ and adjacency given by these edge sets.

	\begin{lemma}\label{lem:orbitcomponents}
The group $\Gamma$ permutes the orbit components, each of which is connected and regular, and is countable under \eqref{eq:countable-orbit-degree}.
The graph $Q$ is connected, and adjacent components lie in different original vertex orbits.
Every small automorphism stabilizes some orbit component setwise and acts nontrivially on it.
	\end{lemma}
\begin{proof}
		The setwise stabilizer of an orbit component acts transitively on its
		vertices. Under \eqref{eq:countable-orbit-degree} its balls, and hence
		the component, are countable. The assertions about $Q$ follow from
		connectedness and the definition of induced components. Finally, a
		witness edge $x\varphi(x)$ lies within one orbit component $H$, forcing
		$\varphi(H)=H$ and a nontrivial restriction.
	\end{proof}

It therefore suffices to arrange
	\begin{equation}\label{eq:componentcriterion}
		\varphi\in\Aut(G,c),\quad\varphi(H)=H
		\quad\Longrightarrow\quad \varphi|_{V(H)}=\id.
	\end{equation}

	\begin{lemma}[Two adjacent anchors]\label{lem:anchors}
Suppose all orbit components belong to $\cH$ and there are at least two.
Choose adjacent components $C,D$, an edge $ab$ between them, and roots $r_H$ with $r_C=a,r_D=b$.
There is a partial list coloring of $E_{\mathrm{int}}\cup E_G(C,D)$ distinguishing every $H$ relative to $r_H$, such that every preserving automorphism of any extension fixes $C,D$ pointwise.
	\end{lemma}
	\begin{proof}
Color $C,D$ relative to their roots using Lemma~\ref{lem:freeedge}.
For each other $H$ with $X(H)=X(C)$, apply Lemma~\ref{lem:freeedge} to choose a rooted distinguishing coloring with no color-preserving isomorphism to the colored graph $C$; do the same with $D$ when $X(H)=X(D)$.
Since $X(C)\ne X(D)$, each $H$ must avoid at most one fixed colored graph.
Color all other components relative to their roots using Lemma~\ref{lem:freeedge}.
Every automorphism maps an orbit component only to one in the same original vertex orbit.
Thus every color preserver fixes $C,D$ setwise.
Give $ab$ a color $p\in L(ab)$ and avoid $p$ on the other edges between them.
This fixes $a,b$ individually, and the internal rooted colorings fix both components pointwise.
	\end{proof}

	\subsection{The 2-list construction}\label{sec:smallproof}
	\begin{proof}[Proof of Theorem~\ref{thm:smalltwo}]
The cases $C_4,K_{3,3}$ follow from Proposition~\ref{prop:c4k33}.
All orbit components are countable and regular by Lemma~\ref{lem:orbitcomponents}.
Fix 2-element edge lists.

\smallskip\noindent\emph{Case 1: some orbit component is outside $\cH$.}
Let $S$ be the union of these components; it is nonempty and $\Gamma$-invariant.
Color them internally by Lemma~\ref{lem:regularsmall}.
For every remaining component $H$, distance to $S$ is constant on $H$, say $d_H>0$, since its vertices belong to one original orbit.
Let $B_H$ be the edges from $H$ to vertices at distance $d_H-1$ from $S$.
Choose $e_H=r_Hu_H\in B_H$, color $H$ distinguishing relative to $r_H$, and arrange
		\begin{equation}\label{eq:uniqueedge}
			c(e_H)=p_H\in L(e_H),\qquad
			c(e)\ne p_H\quad(e\in B_H\setminus\{e_H\}).
		\end{equation}
The sets $B_H$ are pairwise disjoint and contain no internal edge: the endpoint farther from $S$ determines $H$.
Color all remaining edges arbitrarily.
If a preserving automorphism stabilizes $H$, it preserves $B_H$ and fixes its unique $p_H$-edge, hence $r_H$ and then all of $H$.
The other components are already internally distinguished.
Thus \eqref{eq:componentcriterion} holds.

\smallskip\noindent\emph{Case 2: all orbit components belong to $\cH$.}
There are at least two, by the excluded graphs.
Choose adjacent $C,D$ and $ab\in E_G(C,D)$, and put $h(H)=d_Q(H,\{C,D\})$.
For every $H\ne C,D$, choose an edge $e_H=r_Hu_H$ to a component at level $h(H)-1$.
Apply Lemma~\ref{lem:anchors} with these roots and $r_C=a,r_D=b$.
For each remaining $H$ set
		\[
		B_H=\bigcup_{J:\,h(J)=h(H)-1}E_G(H,J)
		\]
and color it as in \eqref{eq:uniqueedge}.
These sets are pairwise disjoint and avoid both internal edges and $E_G(C,D)$; color the remaining edges arbitrarily.
The anchors first force every preserving automorphism to fix $C,D$ pointwise, so it preserves $h$.
If it stabilizes another $H$, it preserves $B_H$, fixes $r_H$ by its unique $p_H$-edge, and fixes $H$ pointwise by its internal coloring.
Again \eqref{eq:componentcriterion} holds.

In either case Lemma~\ref{lem:orbitcomponents} excludes every small color-preserving automorphism.
All assignments use disjoint edge sets, so no recoloring or limiting argument is needed.
	\end{proof}
To pass from connected graphs to arbitrary graphs, we determine the
small edge parameters of a graph from those of its components.
	\begin{lemma}[Connected components]\label{lem:smallcomponents}
If the small edge parameters of the components of a nonempty graph are finite and uniformly bounded, each of $D'_s(G)$ and $D'_{\ell,s}(G)$ is the maximum of the corresponding component parameters.
	\end{lemma}
	\begin{proof}
		Color components independently. A witness edge $x\varphi(x)$ forces
		its connected component to be preserved. Conversely, a small automorphism of a component
		extends by the identity. For the list lower bound, extend a bad list
		assignment arbitrarily to the other components.
	\end{proof}

	\begin{corollary}[Exact small edge parameter]\label{cor:smalledge}
For every nonempty graph $G$ satisfying \eqref{eq:countable-orbit-degree} and having no $K_2$ component,
		\[
		D'_{\ell,s}(G)=D'_s(G)=
		\begin{cases}
			1,&\text{if $G$ has no small automorphism},\\
			3,&\text{if a component is $K_3,K_4,K_5$ or $C_5$},\\
			2,&\text{otherwise}.
		\end{cases}
		\]
	\end{corollary}

\begin{proof}[Proof of Corollary~\ref{cor:smalledge}]
For a component $H$ and $v\in H$, one has $\Aut(G)v\cap V(H)=\Aut(H)v$, so \eqref{eq:countable-orbit-degree} passes to components.
The four exceptional components have both parameters three; every other component has both at most two by Theorem~\ref{thm:smalltwo}, and both equal one exactly when it has no small automorphism.
Apply Lemma~\ref{lem:smallcomponents}.
	\end{proof}

For example, $G_\kappa=K_2\vee\overline{K_\kappa}$, for any infinite cardinal $\kappa$, satisfies \eqref{eq:countable-orbit-degree}, although its two central vertices have degree $\kappa$.
Their transposition is small, so $D'_s(G_\kappa)=D'_{\ell,s}(G_\kappa)=2$.

\section{Further questions}\label{add:sec:questions}
The remaining general questions for connected locally finite graphs, when the ordinary parameters are finite, are
	\[
	D_\ell(G)=D(G),\qquad D'_\ell(G)=D'(G),\qquad
	D_{\ell,s}(G)=D_s(G).
	\]
Theorem~\ref{thm:countable-edge-motion} settles the edge-list infinite-motion conjecture for all countable connected graphs, including those with vertices of infinite degree.
In particular, the middle equality holds whenever every nonidentity automorphism moves infinitely many edges.
This does not settle that equality for graphs with nonidentity automorphisms of finite edge support.
The small edge equality is settled by Corollary~\ref{cor:smalledge} whenever edge distinction is possible, namely in the absence of a $K_2$ component.

Compact transfer controls a rooted stabilizer; the difficulty is to exclude automorphisms moving the root without increasing the list size.
The compact-kernel theorem does this on disjoint invariant coordinate sets.
It remains to determine whether stabilizer transfer extends through an arbitrary compact normal subgroup, whose actions on different finite orbits may be coupled.

\section{AI Declaration}
We used AI to assist with the exploration of proof strategies and the preparation of this manuscript.
The authors take full responsibility for the mathematical content and final text.

	
\end{document}